# From harmonic mappings to Ricci flow due to the Bochner technique

SERGEY STEPANOV[1, 2], IRINA Aleksandrova[1], IRINA TSYGANOK [1]

[1, 2] Department of Data Analysis and Financial Technologies,
Finance University under the Government of Russian Federation,
49-55, Leningradsky Prospect, 125468 Moscow, Russia,
e-mail: *i.i.tsyganok@mail.ru*

[2] Department of Mathematics, Russian Institute for Scientific
and Technical Information of the Russian Academy of Sciences,
20, Usievicha street, 125190 Moscow, Russia,
e-mail: *s.e.stepanov@mail.ru*

**Abstract.** The present paper is devoted to the study a global aspect of the geometry of harmonic mappings and, in particular, infinitesimal harmonic transformations, and represents the application of our results to the theory of Ricci solutions and the Ricci flow. These results will be obtained using the methods of Geometric analysis and, in particular, due to theorems of Yau, Li and Schoen on the connections between the geometry of a complete smooth manifold and the global behavior of its subharmonic functions.

53C20; 53C43; 53C44

## 1. INTRODUCTION

It is well known that there exists a connection between *harmonic maps* and *Ricci flows* due to the "de Turck trick" that modifies the Ricci flow into a nonlinear parabolic equation (see [1, pp. 113-118]). On the other hand, it is well known that self-similar solutions of the Ricci flow are *Ricci solutions* (see [2, p. 22]). At the same time, the vector field that makes a Riemannian metric into a Ricci solution is an *infinitesimal harmonic transformation* (see for the proof [3]; [4]). In this theory we see the unity of geometry and analysis.

In accordance with the above, the main purpose of the present paper is to make an observation of a function-theoretic nature in global differential geometry of harmonic mappings (see [5] - [8]), infinitesimal harmonic transformations (see [9]; [10]), Ricci solutions (see, for example, [2]; [11]) and the Ricci flow (see, for example, [1]). To implement this, our paper is organized as follows: In the second section of the paper, we give a brief survey of basic facts of the geometry "in the large" of harmonic mappings between Riemannian manifolds. In addition, we prove that the classical theorems on harmonic mappings are consequences of well-known assertions on subharmonic functions. Results of the third section of our paper with the title "Infinitesimal harmonic transformations" are obtained as analogs of results of the second section of the paper. In turn, the results of the fourth section which has the title "Ricci solutions" are applications of the results of the third section of our paper. In conclusion, we consider the evolution equations for the scalar curvature and the Ricci tensor under the Ricci flow.

The *Bochner technique* and its generalized version will help us to relate these various research topics. We must recall here that the classic Bochner technique is an analytical method to obtain vanishing theorems for some topological or geometrical invariants on a compact (without boundary) Riemannian manifold, under some curvature assumption (see [12] - [14] and [15, Chapter 9]). The proofs of such theorems reduce to applying the *Bochner maximum principle* and the *theorem of Green* (see [12, p. 30-31]). In the present paper we will also use a generalized version of the Bochner technique (see, for example, [16]) and, in particular, we will use the *maximum principle of Hopf* (see [17]), Yau, Li and Schoen results on the connections between the geometry of a complete smooth manifold and the global behavior of its subharmonic functions (see, for example, [18]; [19]).

Theorems and corollaries of this paper complement our results from [10]; [20] - [23] and the results of other authors from [11]; [13]; [24]; [25, p. 57]; [27].

**Acknowledgments.** Our work was supported by Russian Foundation for Basis Research (projects No. 16-01-00756).

## 2. HARMONIC MAPPINGS FROM THE BOCHNER TECHNIQUE POINT OF VIEW

In this section we shall give brief survey of basic facts of the differential geometry of harmonic mappings between Riemannian manifolds (see the monographs [7] and [8]).

The subject of harmonic maps is vast and has found many applications, and it would require a very long reader to cover all aspects, even superficially (see [8, p. 417]). We have made a choice; in particular, highlighting the key question of nonexistence of harmonic maps between given complete Riemannian manifolds. We will survey some of the main method of Geometric analysis for answering this question (see, for example, [15, Chapter 9]; [16]). In particular, we will prove new versions of two well-known vanishing theorems of harmonic mappings from [5] and [6]. We note here that the question of proving the existence of harmonic mappings was studied in [8].

Let $(M,g)$ and $(\overline{M},\overline{g})$ denote complete Riemannian manifolds of dimensions $n$ and $\overline{n}$, respectively. We call the *energy density* of a smooth map $f:(M,g)\to(\overline{M},\overline{g})$ the nonnegative scalar function $e(f):M\to\mathbb{R}$ such that $e(f)=\frac{1}{2}\|f_*\|^2$ where $\|f_*\|^2$ denotes the squared norm of the differential $f_*$, with respect to the induced metric $\tilde{g}$ on the vector bundle $T^*M\otimes f^*T\overline{M}$ by $g$ and $\overline{g}$ (see also [5]).

It is well known that $f:(M,g)\to(\overline{M},\overline{g})$ is a *harmonic mapping* if and only if it satisfies the *Euler-Lagrange equation*

$$trace_g(\widetilde{\nabla}f_*)=0 \tag{2.1}$$

where $\widetilde{\nabla}=\nabla\oplus\overline{\nabla}$ is the canonical connection in the vector bundle $T^*M\otimes f^*T\overline{M}$ (see [1, p. 117]; [5]; [8]). Moreover, if $f$ is a harmonic mapping then a standard calculation yields (see also [5; p. 123]):

$$\Delta e(f)=\|\widetilde{\nabla}f_*\|^2+Q(f) \tag{2.2}$$

for the *Laplace–Beltrami operator* $\Delta:=div\,\nabla$ and the scalar function

$$Q(f)=g(Ric,f^*\overline{g})-trace_g(trace_g(f^*\overline{R})) \tag{2.3}$$

where $\overline{R}$ is the Riemannian curvature tensor of $(\overline{M},\overline{g})$ and $Ric$ is the Ricci tensor of $(M,g)$.

We showed in [27] the condition when $Q(f)$ is a quasi-positive scalar function on a connected open domain $U\subset M$, i.e. $Q(f)$ is non-negative everywhere in $U$ and $Q(f)$ is positive in at least one point of $U$. Namely, if the sectional curvature $\overline{sec}$ of

$(\overline{M}, \overline{g})$ is nonnegative at an arbitrary point of $f(U) \subset \overline{M}$ and $Ric \geq f^*\overline{Ric}$ at each point of $U$ for the Ricci tensors $\overline{Ric}$ of $(\overline{M}, \overline{g})$, then $Q(f) \geq 0$. Furthermore, if there is at least one point of $U$ at which $Ric > f^*\overline{Ric}$, then $Q(f) > 0$ at this point. In this case, $Q(f)$ is also quasi-positive scalar function defined in $U$. As a result, the energy density function $e(f)$ satisfies the inequality $\Delta e(f) \geq 0$ at each point of $U$, by (2.2). Therefore, $e(f)$ is a *subharmonic function*. The assumption that the energy density function $e(f)$ attains a local maximum value at some point $x \in U$ then implies $e(f)$ is a constant $C$ in $U$, by the *Hopf's maximum principle* (see [12, Theorem 2.1]; [17, Theorem 1]). If $C > 0$, then $grad f$ is nowhere zero. Now, at a point where the $Q(f)$ is positive, the left side of (2.2) is zero while the right side is positive. This contradiction shows that $C = 0$ and hence $f$ is constant in $U$. Thus we have proved the following lemma.

**Lemma 2.1**. *Let $f : (M, g) \to (\overline{M}, \overline{g})$ be a harmonic mapping such that its energy density $e(f)$ has local maximum at some point $x$ in a connected open domain $U \subset M$. If, in addition, the sectional curvature $\overline{sec}$ of $(\overline{M}, \overline{g})$ is nonnegative at an arbitrary point of $f(U) \subset \overline{M}$ and $(M, g)$ has the Ricci tensor $Ric$ such that $Ric \geq f^*\overline{Ric}$ at each point of $U$, and there is at least one point of $U$ at which $Ric > f^*\overline{Ric}$, then f is constant in the domain $U$.*

Using our Lemma 2.1 and the Bochner maximum principle (see [12, Theorem 2.2]), we obtain our theorem from [27] that is a corollary of our Lemma 2.1.

**Theorem 2.2**. *Let $f : (M, g) \to (\overline{M}, \overline{g})$ be a harmonic mapping between Riemannian manifolds $(M, g)$ and $(\overline{M}, \overline{g})$. Assume that the sectional curvature $\overline{sec}$ of the second manifold $(\overline{M}, \overline{g})$ is nonnegative in every point of $f(M)$ and the first manifold $(M, g)$ is a compact manifold with the Ricci tensor $Ric \geq f^*\overline{Ric}$. Then f is totally geodesic mapping with constant energy density $e(f)$. Furthermore, if there is at least one point of M in which $Ric > f^*\overline{Ric}$, then f is a constant mapping.*

**Remark 2.1.** We recall here that in the well known paper of Eells and Sampson [5] was proved the celebrated vanishing theorem on harmonic maps which state the following: If $f:(M,g)\to(\overline{M},\overline{g})$ is any harmonic mapping between a compact Riemannian manifold (M, g) with the Ricci tensor $Ric \geq 0$ and a Riemannian manifold $(\overline{M},\overline{g})$ with the sectional curvature $\overline{sec} \leq 0$ then f is *totally geodesic* and has constant the *energy density* $e(f)$. Furthermore, if there is at least one point of M at which its Ricci curvature $Ric > 0$, then every harmonic map $f:(M,g)\to(\overline{M},\overline{g})$ is constant.

Li and Schoen have proved in [19] that there is no non-constant, non-negative $L^p$-integrable $(0<p<\infty)$ subharmonic function $u$ on any complete Riemannian manifold $(M,g)$ with nonnegative Ricci tensor. In other word, if we suppose that $Ric \geq 0$ and $\int_M \|u\|^p dVol_g < \infty$ for a complete Riemannian manifold $(M,g)$, then $u=C$ for some constant $C$. In this case, we have $C^p \int_M dVol_g < \infty$. If $C>0$, $u$ is nowhere zero and the volume of $(M,g)$ is finite. At the same time, we know from [18] that every complete non-compact Riemannian manifold $(M,g)$ with non-negative Ricci tensor has infinite volume. This contradiction shows $C=0$ and hence $u \equiv 0$. Therefore, we can formulate the following

**Lemma 2.3**. *Let $(M,g)$ be a complete non-compact Riemannian manifold with nonnegative Ricci tensor, then there is no nonzero non-negative $L^p$-integrable $(0<p<\infty)$ subharmonic function on $(M,g)$.*

Let $(M,g)$ be a complete non-compact Riemannian manifold. Given a smooth map $f:(M,g)\to(\overline{M},\overline{g})$ we definite its *energy* as (see [5]; [6]):

$$E(f) = \int_M e(f) \, dVol_g \, .$$

The energy $E(f)$ can be infinite or finite for a smooth map $f:(M,g)\to(\overline{M},\overline{g})$ of a complete Riemannian manifold $(M,g)$. In particular, $E(f) < +\infty$ for a compact Riemannian manifold $(M,g)$ (see [6]).

At the same time, we can formulate an alternative theorem for harmonic maps from complete Riemannian manifolds to Riemannian manifolds with nonnegative sectional curvature (see [6] and [7, p. 25]). In this case, we can also use Lemma 2.1 on nonzero non-negative $L^p$-integrable $(0 < p < \infty)$ subharmonic functions on a complete Riemannian manifold with nonnegative Ricci curvature for the proof of this theorem. Namely, the following theorem holds.

**Theorem 2.4.** *Let $f : (M, g) \to (\overline{M}, \overline{g})$ be a harmonic mapping with finite energy. If $(M, g)$ is a complete non-compact Riemannian manifold with Ricci tensor $Ric \geq f^* \overline{Ric}$, where $\overline{Ric}$ denotes the Ricci tensor of $(\overline{M}, \overline{g})$, and $(\overline{M}, \overline{g})$ is a Riemannian manifold with nonnegative sectional curvature $\overline{sec}$ in every point of $f(M)$, then f is a constant map.*

**Remark 2.2.** We recall here that Yau and Schoen showed the following well-known vanishing theorem (see [6]): A harmonic map of finite energy $E(f)$ from a complete non-compact manifold $(M, g)$ with the Ricci tensor $Ric \geq 0$ to a compact manifold $(\overline{M}, \overline{g})$ with the sectional curvature $\overline{sec} \leq 0$ is homotopic to a constant map.

## 2. INFINITESIMAL HARMONIC TRANSFORMATIONS FROM THE FUNCTION-THEORETIC POINT OF VIEW

The main results of this section are obtained as analogs of results of the second section of the present paper.

A vector field $\xi$ on a complete Riemannian manifold $(M, g)$ is called an *infinitesimal harmonic transformation* (see [9]) if $\xi$ generates a flow which is a local one-parameter group of harmonic transformations $(t, x) \in \mathbb{R} \times M \to \varphi_t(x) \in M$ (in the other words, local harmonic diffeomorphisms). Analytic characteristic of such vector field has the form $trace_g(L_\xi \nabla) = 0$ for the Lie derivative $L_\xi$ in the direction of $\xi$ (see [9]; [10]). This formula is an analogue of the formula (2.1). In addition, we have proved in [10] that a vector field $\xi$ is an infinitesimal harmonic transformation if and only if $\tilde{\Delta}\theta = 2\,Ric(\xi,\cdot)$ for the 1-form $\theta$ corresponding to $\xi$ under the duality defined by the metric $g$ and the *Hodge-de Rham Laplacian* $\tilde{\Delta}$ (see [14, p. 158]).

In accordance with the theory of *harmonic maps* (see [5]) we define the *energy density* of the flow on $(M, g)$ generated by an infinitesimal harmonic transformation $\xi$ as the scalar function $e(\xi) = \frac{1}{2}\|\xi\|^2$ where $\|\xi\|^2 = g(\xi, \xi)$. Then the Laplace–Beltrami operator $\Delta e(\xi)$ for the energy density $e(\xi)$ of an infinitesimal harmonic transformation $\xi$ has the form (see [4]; [21]; [22])

$$\Delta e(\xi) = \|\nabla \xi\|^2 - Ric(\xi, \xi). \tag{3.1}$$

The formula (3.1) is an analogue of the formula (2.2) for the energy density $e(f)$ of a harmonic map $f$. In this case, the following theorem is true.

**Theorem 3.1**. *Let $(M, g)$ be a Riemannian manifold and $U \subset M$ be a connected, open domain. If the energy density of the flow $e(\xi)$ generated by an infinitesimal harmonic transformation $\xi$ has a local maximum in some point of $U$ and the Ricci curvature of $(M, g)$ is quasi-negative in $U$, then $\xi \equiv 0$ everywhere in $U$.*

**Proof.** Let the Ricci curvature of $(M, g)$ is quasi-negative everywhere in a connected, open domain $U \subset M$, then the energy density function $e(\xi)$ satisfies the inequality $e(\xi) \geq 0$, by (3.8). This means that $e(\xi)$ is a subharmonic function. Suppose now that the energy density function $e(\xi)$ attains a local maximum value in a point $x \in U$, then $e(\xi)$ is a constant $C$ in $U$, by Hopf's maximum principle (see [12, Theorem 2.1]; [17, Theorem 1]). If $C > 0$, $\xi$ is nowhere zero. Now, at a point where the $Ric$ is negative, the left side of (3.1) is zero while the right side is positive. This contradiction shows $C = 0$ and hence $\xi \equiv 0$ everywhere in the domain $U$.

**Remark 3.1.** Our Theorem 3.1 is a direct generalization of Theorem 4.3 presented in Kobayashi's monograph on transformation groups (see [25, p. 57]) and Wu's proposition on a Killing vector field whose length achieves a local maximum (see [13]).

As an analogue of Theorem 2.2, we can formulate the following theorem, which can be proved using the *Bochner maximum principal* (see [12, Theorem 2.2]).

**Theorem 3.2**. *A compact Riemannian manifold $(M, g)$ with quasi-negative Ricci curvature doesn't have nonzero infinitesimal harmonic transformation.*

Next, we recall that the *kinetic energy* $E(\xi)$ of the flow on $(M, g)$ generated by a vector field $\xi$ is determined in accordance with [28, p. 2] by the following equation:

$$E(\xi) = \int_M e(\xi) \, dVol_g.$$

**Remark 3.2.** The definition is consistent with the theory of harmonic mappings in the case of an infinitesimal harmonic transformation. Moreover, the energy $E(\xi)$ can be infinite and finite. For example, $E(\xi) < +\infty$ for a smooth complete vector field $\xi$ on a compact Riemannian manifold $(M, g)$.

As an analogue of Theorem 2.4, we formulate the following

**Theorem 3.3**. *Let $(M, g)$ be a complete Riemannian manifold $(M, g)$ with the nonpositive Ricci curvature. Then every infinitesimal harmonic transformation with finite kinetic energy is parallel. If, in addition, the volume of $(M, g)$ is infinite or the Ricci curvature is negative at some point of $M$, then each infinitesimal harmonic transformation is identically zero everywhere on $(M, g)$.*

**Proof.** For the proof we use the well known *second Kato inequality* (see [14, p. 380])

$$-\|\xi\| \Delta \|\xi\| \leq g(\overline{\Delta}\theta, \theta)$$

where $\overline{\Delta} := -trace_g \nabla \circ \nabla$ is the *rough Laplacian* and $\theta$ is the 1-form corresponding to $\xi$ under the duality defined by the metric $g$. It is well known that the rough Laplacian satisfies the *Weitzenböck formula* (see [25, p. 44]; [14, p. 378])

$$\overline{\Delta}\theta = \widetilde{\Delta}\theta - S\xi$$

where $S$ is the *Ricci operator* defined by $g(SX, Y) = Ric(X, Y)$ for any tangent vector fields $X$ and $Y$. Therefore, the second Kato inequality can be rewritten in the form

$$2\sqrt{e(\xi)} \Delta \sqrt{e(\xi)} \geq -g(\widetilde{\Delta}\theta, \theta) + Ric(\xi, \xi). \qquad (3.2)$$

where $\|\xi\| = \sqrt{2e(\xi)}$. On the other hand, we have proved in [10] and [21] that $\xi$ is an infinitesimal harmonic transformation on $(M, g)$ if and only if $\widetilde{\Delta}\xi = 2S\xi$. Therefore, we obtain from (3.2) the following equation

$$\sqrt{e(\xi)}\,\Delta\sqrt{e(\xi)} = -\frac{1}{2}Ric(\xi,\xi). \tag{3.3}$$

Let us require the Ricci tensor $Ric$ to be nonpositive. In [18, p. 664] and [29] it was shown that every nonnegative smooth function $u$ defined on a complete Riemannian manifold $(M,g)$ and satisfying the conditions $u\Delta u \geq 0$ and $\int_M u^p\,dVol_g < +\infty$ for all $p \neq 1$, must be constant. In particular, if the volume of $(M,g)$ is infinite and $u = const$, then $u = 0$. Therefore, if the smooth manifold $(M,g)$ is complete and

$$E(\xi) = \int_M e(\xi)\,dVol_g < +\infty, \tag{3.4}$$

then the function $\sqrt{e(\xi)}$ is constant. At the same time, we obtain from (3.4) that the volume of $(M,g)$ is finite. Thus it follows from (3.2) that $\nabla\xi = 0$. On the other hand, if we suppose that $Ric_x < 0$ in a point $x \in M$, then this inequality contradicts the equation (3.1). The proof is complete.

**Remark 3.4.** If $\xi$ is nonzero vector field such that $Ric(\xi,\xi) \leq 0$, then (3.3) implies the inequality $\Delta\sqrt{e(\xi)} \geq 0$. This means that $\sqrt{e(\xi)}$ is a subharmonic function. Therefore, the result given above is an analogue of the result which we showed in Theorem 2.3.

The following corollary is valid.

**Corollary 3.1.** *Let $(M,g)$ be a connected complete noncompact Riemannian manifold of dimension $n \geq 2$ with irreducible holonomy group $\mathrm{Hol}(g)$ and nonpositive Ricci curvature. Then each infinitesimal harmonic transformation on $(M,g)$ with finite kinetic energy is identically zero everywhere on $(M,g)$.*

**Proof.** By Theorem 3.3, an infinitesimal harmonic transformation $\xi$ with $E(\xi) < +\infty$ on a connected complete noncompact Riemannian manifold $(M,g)$ with nonpositive Ricci curvature is parallel. Under the assumption that holonomy group $\mathrm{Hol}(g)$ is irreducible, this relation means that $\xi \equiv 0$.

An example of complete smooth manifold $(M,g)$ with nonpositive Ricci curvature is the well know *Cartan-Hadamard manifold*, that is a simply connected complete Riemannian manifold of nonpositive sectional curvature. The following assertion is valid.

**Corollary 3.2.** *Let $(M, g)$ be a Cartan-Hadamard manifold of dimensional $n \geq 2$ with infinite volume, then each infinitesimal harmonic transformation on $(M, g)$ with finite kinetic energy is identically zero everywhere on $(M, g)$.*

**Remark 3.5.** Other properties of infinitesimal harmonic transformations can be found in our papers [4]; [10]; [21] - [25]. In particular, we proved in [10] that the set of all infinitesimal harmonic transformations on a compact Riemannian manifold $(M, g)$ is a finite-dimensional vector space over $\mathbb{R}$. Moreover, the Lie algebra $i(M)$ of infinitesimal isometric transformation is a subspace of this vector space (see [10]). We recall here that an infinitesimal isometric transformation (infinitesimal isometry) or a Killing vector field $X$ on $(M, g)$ is defined by the well known equation $L_X g = 0$. It is well known that $\dim i(M) = 1/2 \, n(n+1)$ on an $n$-dimensional Riemannian manifold $(M, g)$ of constant curvature (see, for example, [25, pp. 46-47]). Therefore, the dimension of the vector space of infinitesimal harmonic transformations on an $n$-dimensional Riemannian manifold $(M, g)$ of constant curvature is at least $1/2 \, n(n+1)$. It should be noted that this proposition is a local result.

The following statement on infinitesimal isometric transformations is a well known oldest result (see, for example, [12, p. 57] and [25, p. 44]): If $\xi$ is an infinitesimal isometry, it satisfies the following differential equations:

$$\widetilde{\Delta} \theta = 2 \, Ric(\xi, \cdot); \tag{3.5}$$

$$div \, \xi = 0 \tag{3.6}$$

for the 1-form $\theta$ corresponding to $\xi$ under the duality defined by the metric $g$. Conversely, if $(M, g)$ is compact and $\xi$ satisfies (3.5) and (3.6), then $\xi$ is an infinitesimal isometry.

The equation (3.6) is more rigid than required for the statement above. In turn, we can formulate and prove an alternative version of this statement.

**Theorem 3.3.** *Let $(M, g)$ be a Riemannian manifold and $\xi$ a vector field on $(M, g)$. If $\xi$ is an infinitesimal isometry, then it satisfies the following conditions: $\xi$ is an infinitesimal harmonic transformation and*

$$L_\xi \, div\, \xi \geq 0. \tag{3.7}$$

*Conversely, if $(M, g)$ is compact and $\xi$ satisfies the above two conditions, then $\xi$ is an infinitesimal isometry.*

**Proof.** We have already shown that if $\xi$ is an infinitesimal harmonic transformation, it satisfies (3.5). If, in addition, $\xi$ is an infinitesimal isometry, then the equation $div\,\xi = 0$ holds. This implies (3.7). To prove the converse, we may assume that $M$ is compact and orientable. If $M$ is not orientable, consider its orientable double covering. Let's consider the vector field $X = (div\,\xi)\xi$ for an arbitrary infinitesimal harmonic transformation $\xi$ on $(M, g)$. The divergence of $X$ has the form

$$div\, X = L_\xi (div\,\xi) + (div\,\xi)^2. \tag{3.8}$$

If we apply classic Green's theorem $\int_M div\, X\, dVol_g = 0$ to $X = (div\,\xi)\xi$, then we obtain the integral formula

$$\int_M \left( L_\xi (div\,\xi) + (div\,\xi)^2 \right) dVol_g = 0 \tag{3.9}$$

for the canonical measure $dVol_g$ which is associated to the metric $g$. If the inequality $L_\xi(div\,\xi) \geq 0$ holds anywhere on $(M, g)$, then from (3.9) we conclude that $div\,\xi = 0$. Next, for complete the proof we can refer to Theorem 3.4.

*Remark 3.5.* We can define the divergence of a vector field using a dynamics approach. Namely, the divergence of each smooth vector field $\xi$ on $(M, g)$ is a scalar function defined by (see, for example, [25, p. 6] and [30, p. 195])

$$(div\,\xi)dVol_g = L_\xi(dVol_g). \tag{3.10}$$

Specifically, the formula (3.10) checking how the volume form changes along the flow of the vector field. Due to (3.10), the function $div\,\xi$ is called in [30, p. 195] the *logarithmic rate of volumetric expansion* along the flow generated by the vector field $\xi$. Therefore, $L_\xi \, div\, \xi$ measures the *acceleration of volumetric expansion*, i.e. the acceleration of change of the volume element $dVol_g$ along trajectories of the flow with the velocity vector $\xi$. In particular, the condition $L_\xi \, div\, \xi \geq 0$ means that $dVol_g$ is a nondecreasing function along trajectories of this flow.

**Remark 3.6.** The equality $L_\xi g = 0$ also implies the invariance condition for the Ricci tensor $L_\xi Ric = 0$. In General Relativity, there are investigations (see, for example, [31] and [32]), where "weakened condition" of the form $trace_g(L_\xi Ric) = 0$ is studied instead of the condition $L_\xi Ric = 0$.

It is of interest to note, that, in the case of compact $(M, g)$, the addition of conditions $L_\xi Ric \geq 0 \, (\leq 0)$ to equation (3.1) implies that infinitesimal harmonic transformation $\xi$ will actually be an infinitesimal isometry. Namely, the following theorem holds.

**Theorem 3.4.** *Let $(M, g)$ be a Riemannian manifold and $\xi$ a vector field on $(M, g)$. If $\xi$ is an infinitesimal isometry, it satisfies the following conditions: $\xi$ is an infinitesimal harmonic transformation and*

$$trace_g(L_\xi Ric) \geq 0 \, (\leq 0). \qquad (3.11)$$

*Conversely, if $(M, g)$ is compact and $\xi$ satisfies the above two conditions, then $\xi$ is an infinitesimal isometry.*

**Proof.** We have already shown that if $\xi$ is an infinitesimal harmonic transformation, it satisfies (3.5) that equals to the first condition of our theorem. If, in addition, $\xi$ is an infinitesimal isometry, then the equation $L_\xi Ric = 0$ holds. To prove the converse, we may assume that $M$ is compact. If $\xi$ is an infinitesimal harmonic transformation then we have the differential equation (see [33])

$$\Delta\, div\, \xi = trace_g(L_\xi Ric).$$

Using the *Bochner maximum principal* (see [12, Theorem 2.2]), we can conclude that $div\, \xi = const$. On the other hand, we have $\int_M div\, \xi \, dVol_g = 0$. Hence, $div\, \xi = 0$, showing that $\xi$ is an infinitesimal transformation (see Theorem 3.4). Theorem 3.6 is proved.

**Remark 3.7.** Theorem 3.4 is our second alternative version of the classic theorem on infinitesimal isometric transformations (see, for example, [12, p. 57] and [25, p. 44]).

# 4. RICCI SOLITONS FROM INFINITESIMAL HARMONIC TRANSFORMATIONS POINT OF VIEW

The main results of this section are applications to the Ricci soliton theory of the results of the third section of our paper.

Let $g$ be a fixed Riemannian metric on a smooth manifold $M$. Consider the one-parameter family of diffeomorphisms $(t, x) \in \mathbb{R} \times M \to \varphi_t(x) \in M$ that is generated by a smooth vector field $\xi$ on $M$. The evolutive metric $g(t) = \sigma(t)\varphi_t^*(x)g(0)$ for a positive scalar $\sigma(t)$ such that $\sigma(0) = 1$ and $g(0) = g$ is a *Ricci soliton* iff the metric $g$ is a solution of the nonlinear stationary PDF

$$-2\,Ric = L_\xi\, g + 2\lambda g \qquad (4.1)$$

where $Ric$ is the Ricci tensor of $g$, $L_\xi\, g$ is the *Lie derivative* of $g$ with respect to $\xi$ and $\lambda$ is a constant (see, for example, [2, p. 22]). To simplify the notation, we denote the Ricci soliton in the following way $(g,\xi,\lambda)$. A Ricci soliton is called *steady*, *shrinking* and *expanding* if $\lambda = 0$, $\lambda < 0$ and $\lambda > 0$, respectively. In addition, a Ricci soliton is called *Einstein* if $L_\xi\, g = 0$, and it is called *trivial* if $\xi \equiv 0$.

In [3], we have shown that the following theorem is true.

**Theorem 4.1**. *The vector field $\xi$ of an arbitrary Ricci soliton $(g,\xi,\lambda)$ on a smooth manifold $M$ is an infinitesimal harmonic transformation on the Riemannian manifold $(M, g)$.*

Therefore, from Theorem 3.1 and Theorem 4.1 we can obtain the following proposition.

**Corollary 4.2**. *Let $(g,\xi,\lambda)$ be a Ricci soliton on a smooth manifold $M$. If the Ricci tensor $Ric$ of $g$ is quasi-negative in a connected, open domain $U \subset M$ and the energy density of the flow generated by the vector field $\xi$ has a local maximum in some point of $U$, then $(g,\xi,\lambda)$ is an expanding Ricci soliton.*

**Proof.** Let $(g,\xi,\lambda)$ be a Ricci soliton on a smooth manifold $M$, then its vector field $\xi$ is an infinitesimal harmonic transformation on the Riemannian manifold $(M, g)$. If, in addition, the Ricci tensor $Ric$ of $g$ is quasi-negative in a connected, open domain $U \subset M$ and the energy density of the flow $e(\xi)$ generated by $\xi$ has a local maximum in some point $x$ of $U$, then $\xi \equiv 0$ in any point $y \in U$, by our Theorem 3.1. In this case, from (4.1) we obtain the equality $Ric_y = -\lambda g_y$ and hence $\lambda > 0$.

**Remark 4.1.** From Corollary 3.1 we can conclude that a Ricci soliton $(g,\xi,\lambda)$ with the quasi-negative Ricci curvature of $g$ is trivial on a compact smooth manifold $M$.

In turn, from Theorem 3.2 we obtain the following

**Corollary 4.3.** *Let $(g,\xi,\lambda)$ be a Ricci soliton on a compact smooth manifold $M$. If the volume element $dVol_g$ (resp. the scalar curvature $s$ of $g$) is a nondecreasing (resp. nonincreasing) function along trajectories of the flow with the velocity vector $\xi$, then $(g,\xi,\lambda)$ is trivial.*

**Proof.** Consider a Ricci soliton $(g,\xi,\lambda)$ on a compact smooth manifold $M$. From (4.1) we obtain $div\,\xi = -(s+n\lambda)$ for the scalar curvature $s = trace_g\,Ric$, then the acceleration of volumetric expansion of the flow generated by the vector field $\xi$ of the Ricci soliton $(g,\xi,\lambda)$ has the form

$$L_\xi(div\,\xi) = -L_\xi s. \tag{4.2}$$

Therefore, the acceleration of change of the volume element $dVol_g$ along trajectories of the flow with the velocity vector $\xi$ equals to $-L_\xi s$. In particular, the condition $L_\xi s \leq 0$ means that $dVol_g$ is a nondecreasing function along trajectories of this flow. On the other hand, if $L_\xi s \leq 0$ then from Theorem 3.5 we obtain that $\xi$ is an infinitesimal isometry, i.e. $L_\xi g = 0$. On the other hand, an arbitrary Ricci soliton $(g,\xi,\lambda)$ on a compact smooth manifold $M$ is a gradient soliton, i.e. $\theta = grad\,u$ for some $u \in C^\infty M$ (see [11]). Therefore, the condition $L_\xi g = 0$ becomes $\nabla\nabla u = 0$ which implies the equation $\Delta u = 0$, i.e. $u$ is a *harmonic function*. In this case, $u = const$, by Bochner maximum principal. As a result, we obtain $\theta = grad\,u = 0$, and hence, the Ricci soliton $(g,\xi,\lambda)$ is trivial.

The following corollary of Theorem 3.6 is proved similarly.

**Corollary 4.4.** *Let $(g,\xi,\lambda)$ be a Ricci soliton on a compact smooth manifold $M$. If $trace_g\, L_\xi Ric \geq 0\ (\leq 0)$ for the Ricci tensor $Ric$ of $g$, then $(g,\xi,\lambda)$ is trivial.*

**Remark 4.2.** It is well known, that every steady and expanding Ricci soliton on a compact smooth manifold $M$ is trivial (see, for example, [11]). On the other hand, the well known problem presented in [11]: are there special conditions in dimensional $n \geq 4$ assuring that a shrinking compact Ricci soliton is trivial? Our two corollaries are answers to this problem.

The proof of the following corollary which is obtained from Theorems 3.3 is not required.

**Corollary 4.5.** *Let $M$ be a connected smooth manifold and $(g,\xi,\lambda)$ be a Ricci soliton with complete Riemannian metric $g$ and nonpositive Ricci curvature on $M$. If the kinetic energy of the flow generated by the vector field $\xi$ is infinite, then $(g,\xi,\lambda)$ is an Einstein Ricci soliton. If, in addition, the volume of $(M,g)$ is infinite or the Ricci curvature is negative at some point, then $(g,\xi,\lambda)$ is a trivial Ricci soliton.*

From (3.8) and (4.2) we conclude that the divergence of the acceleration vector of volumetric expansion $X = (div\,\xi)\xi$ has the form

$$div\, X = -L_\xi s + (s + n\lambda)^2. \qquad (4.3)$$

If, in addition, $(M,g)$ is a complete and oriented Riemannian manifold such that $\|(div\,\xi)\xi\| \in L^1(M,g)$ and $L_\xi s \leq 0$, then from (4.3) we obtain that $div\, X \geq 0$. Then, thanks to a *generalized Green's theorem* (see [34]; [35]), we have $div\, X = 0$. It means that $L_\xi s = 0$ and

$$s = -n\lambda. \qquad (4.4)$$

In this case, from the well-known Schur's identity $\delta Ric = -2^{-1}\nabla s$, we obtain the equation $\delta Ric = 0$. We recall the well-known equality $\delta \widetilde{\Delta} = \widetilde{\Delta}\delta$. Therefore, if we apply the divergence operator $\delta$ to both sides of the equation $\widetilde{\Delta}\theta = 2\,Ric(\xi,\cdot)$, we obtain $trace_g(S \cdot \nabla \xi) = 0$. In this case, from the equation (4.1) and the equality (4.4) we obtain the following $\|Ric\|^2 = n^{-1}s^2$. Then the square of the traceless part of the Ricci tensor is equal to zero, i.e. $\|Ric - n^{-1}s\,g\|^2 = \|Ric\|^2 - n^{-1}s^2 = 0$. Therefore, we conclude that $Ric = n^{-1}s \cdot g$ and $L_\xi g = 0$. This means that $\xi$ is an infinitesimal isometry and $(g,\xi,\lambda)$ is an Einstein soliton. Thus, we proved the following statement.

**Theorem 4.6.** *Let $(g,\xi,\lambda)$ be a Ricci soliton with complete Riemannian metric $g$ on a connected and oriented smooth manifold M such that*

(i) *the volume element $d\,Vol_g$ (resp. the scalar curvature $s$ of $g$) is a nondecreasing (resp. nonincreasing) function along trajectories of the flow with the velocity vector $\xi$;*

(ii) *$\|X\| \in L^1(M,g)$ for the logarithmic rate of volumetric expansion $X = (div\,\xi)\xi$, then the flow generated by $\xi$ consists of isometric transformations and $(g,\xi,\lambda)$ is an Einstein soliton.*

**Remark 4.3.** Our Theorem 4.6 generalizes Theorem 5.4 which is one of the main results of the paper [24].

In particular, for the case $s = $ constant we have the following

**Corollary 4.7.** *Let $(g,\xi,\lambda)$ be a Ricci soliton with complete Riemannian metric $g$ on a connected and oriented smooth manifold M such that $\|X\| \in L^1(M,g)$ for the logarithmic rate of volumetric expansion $X = (div\,\xi)\xi$, then the flow generated by $\xi$ consists of isometric transformations and $(g,\xi,\lambda)$ is an Einstein soliton.*

**Remark 4.4.** The Corollary 4.7 complements the results of the article [26], where the authors studied complete gradient Ricci solitons which have constant scalar curvature. Let $(g,\xi,\lambda)$ be a shrinking Ricci soliton with complete metric $g$ on a connected and oriented smooth manifold $M$, such that $\xi$ satisfies (i) and (ii), then from Theorem 4.2 we obtain the following inequality $Ric = n^{-1}s \cdot g > 0$ for the positive constant $s$. Then it follows from Myers' diameter bound that $(M,g)$ must be compact (see, for example, [15, pp. 251 and 386]). Suppose that $\dim M = 3$, then $(M,g)$ has a constant sec-tional curvature and, in particular, if $M$ is simply connected it be isometric a Euclidian sphere $S^3$. In this case, we have the following

**Corollary 4.8**. *Let $(g,\xi,\lambda)$ be a shrinking Ricci soliton with complete Riemannian metric $g$ on a connected and oriented simply connected three-dimensional smooth manifold M such that the volume element $dVol_g$ is a nondecreasing function along trajectories of the flow with the velocity vector $\xi$ and $\|X\| \in L^1(M,g)$ for the logarithmic rate of volumetric expansion $X = (div\,\xi)\xi$, then $(M,g)$ be isometric a Euclidian sphere $S^3$.*

In conclusion, Corollary 3.1 implies the following assertion.

**Corollary 4.9**. *Let $(g,\xi,\lambda)$ be a Ricci soliton with complete metric $g$, nonpositive Ricci curvature and irreducible holonomy group $\mathrm{Hol}(g)$ on a connected noncompact smooth manifold $M$. If the kinetic energy of the flow generated by the vector field $\xi$ is infinite, then $(g,\xi,\lambda)$ is a trivial soliton.*

## 5. ON THE EVOLUTION EQUATIONS OF THE SCALAR CURVATURE AND THE RICCI TENSOR

Given a 1-parametric family of metrics $g(t)$ on a manifold $M$, defined on a time interval $J \subset \mathbb{R}$, then the Hamilton's *Ricci flow* equation is

$$\frac{\partial}{\partial t}g(t) = -2Ric \tag{5.1}$$

For the Ricci tensor $Ric$ of the metric $g = g(0)$. It is well known that for any $C^\infty$-metric $g$ on a compact manifold $M$, the exists a unique solution $g(t)$, $t \in [0,\varepsilon)$, to the Ricci flow equation for some $\varepsilon > 0$, with $g(0) = g$ (see [36]).

**Remark 5.1.** Compact Ricci solitons are the fixed point of the Ricci flow (5.1) projected from the space of metrics onto its quotient modulo diffeomorphisms and scalings, and often arise as blow-up limits for the Ricci flow on compact manifolds.

Under the Ricci flow, we have an *evolution equation* for the scalar curvature $s$ (see [1, p. 99])

$$\frac{\partial}{\partial t}s = \Delta s + 2\|Ric\|^2. \tag{5.2}$$

If we suppose here that $\frac{\partial}{\partial t}s \leq 0$ for any $t \in J$, then from (5.2) we obtain

$$\Delta s \leq -2\|Ric\|^2 \leq 0. \tag{5.3}$$

For the case of a compact manifold $M$ we obtain $\Delta s = 0$ and $s = constant$ by the "Bochner maximum principal". Then from (5.3) we conclude that $Ric \equiv 0$. In this case, the equation (5.2) can be rewritten in the form $\frac{\partial}{\partial t}g \equiv 0$ and hence, the Ricci flow (5.1) is trivial. Then we can formulate the following obvious proposition.

**Proposition 5.1**. *Let $M$ be a compact smooth manifold with the Hamilton's Ricci flow $\frac{\partial}{\partial t}g(t) = -2Ric$ for a 1-parametric family of metrics g(t) on M, defined on a time interval $J \subset \mathbb{R}$. If $\frac{\partial}{\partial t}s \leq 0$ for any $t \in J$ then the Ricci flow is trivial.*

**Remark 5.1.** Our Corollary 4.3 automatically follows from this proposition.

On the other hand, the evolution equation of the Ricci tensor under the Ricci flow has the form (see [1, p. 112])

$$\frac{\partial}{\partial t}Ric = \Delta_L Ric \tag{5.4}$$

where $\Delta_L$ denotes *the Lichnerowicz Laplacian* (see [1, p. 109]). From (5.4) we obtain

$$trace_g\left(\frac{\partial}{\partial t}Ric\right) = \Delta s. \tag{5.5}$$

If we suppose here that $trace_g\left(\frac{\partial}{\partial t}Ric\right) \geq 0$ or $trace_g\left(\frac{\partial}{\partial t}Ric\right) \leq 0$ for any $t \in J$, then from (5.5) we obtain $\Delta s \leq 0$ or $\Delta s \geq 0$, respectively. For the case of a compact manifold $M$ we obtain $s = constant$ by the "Bochner maximum principal". In this case, the Ricci flow (5.1) is trivial. Then we can formulate the following obvious proposition.

**Proposition 5.2.** *Let M be a compact smooth manifold with the Hamilton's Ricci flow $\frac{\partial}{\partial t} g(t) = -2 Ric$ for a 1-parametric family of metrics g(t) on M, defined on a time interval $J \subset \mathbb{R}$. If $trace_g \left( \frac{\partial}{\partial t} Ric \right) \geq 0$ or $trace_g \left( \frac{\partial}{\partial t} Ric \right) \leq 0$ for any $t \in J$ then the Ricci flow is trivial.*

**Remark 5.2.** Our Corollary 4.4 automatically follows from this proposition.

## REFERENCES

[1] Chow B., Lu P., Ni L., *Hamilton's Ricci Flow*, AMS (2006).

[2] Chow B., Knopf D., *The Ricci flow: An introduction*, American Mathematical Society, Providence (2004).

[3] Stepanov S.E., Shelepova V.N., *A note on Ricci soliton*, Math. Notes 86 (3), 447-450 (2009).

[4] Stepanov S.E., Mikes J., *The spectral theory of the Yano Laplacian with some of its applications*, Ann. Glob. Anal. Geom. 48, 37-46 (2015).

[5] Eells, J., Sampson, J.H., *Harmonic mappings of Riemannian manifolds*, American Journal of Mathematics, 86 (1), 109-160 (1964).

[6] Schoen R., Yau S.T., *Harmonic maps and the topology of stable hypersurfaces and smooth manifolds with nonnegative Ricci curvature*, Comment. Math. Helvetici 61 (1), 333-341 (1976).

[7] Xin Y., *Geometry of harmonic maps*, Springer Science and Business Media, Boston and Berlin (2012).

[8] Hélein F., Wood J., *Harmonic maps*, Handbook of global analysis, 417-491, Elsevier Sci. B.V., Amsterdam (2008).

[9] Nouhaud O., *Transformations infinitésimales harmoniques*, C. R. Acad. Sc. Paris, 274, 573–576 (1972).

[10] Stepanov S.E., Shandra I.G., *Geometry of infinitesimal harmonic transformation*, Ann. Glob. Anal. Geom., 24, 291-297 (2003).

[11] Eminenti M., La Nave G., Mantegazza C., *Ricci solitons: The equation point of view*, Manuscripta Mathematica, 127 (3), 345-367 (2008).

[12] Bochner S., Yano K., *Curvature and Betti numbers*, Princeton, Princeton University Press (1953).

[13] Wu H., *A remark on the Bochner technique in differential geometry*, Proc. Amer. Math. Soc., 78 (3), 403-408 (1980).

[14] Bérard P.H., *From vanishing theorems to estimating theorems: the Bochner technique revisited*, Bulletin of the American Math. Soc., 19 (2), 371-406 (1988).

[15] Peterson P., *Riemannian geometry*, Springer Int. Publ. AG, Switzerland (2016).

[16] Pigola S., Rigoli M., Setti A.G., *Vanishing and Finiteness Results in Geometric Analysis. A Generalization of the Bochner Technique*, Birkhäuser Verlag AG, Berlin (2008).

[17] Calabi E., *An extension of E. Hopf's maximum principle with an application to Riemannian geometry*, Duke Math. J., 25, 45-56 (1957).

[18] Yau S.T., *Some function-theoretic properties of complete Riemannian manifold and their applications to geometry*, Indiana Univ. Math. J., 25 (7), 659-679 (1976).

[19] Li P., Schoen R., *$L^p$ and mean value properties of subharmonic functions on Riemannian manifolds*, Acta Mathematica, 153, 279-301 (1984).

[20] Stepanov S.E., Tsyganok I.I. Mikesh J., *From infinitesimal harmonic transformations to Ricci solitons*, Mathematica Bohemica, 1, 25-36 (2013).

[21] Stepanov S.E., Shandra I.G., *Harmonic diffeomorphisms of smooth manifolds*, St. Petersburg. Math. J., 16 (2), 401–412 (2005).

[22] Stepanov S.E., Tsyganok I.I., *Infinitesimal harmonic transformations and Ricci solitons on complete Riemannian manifolds*, Russian Math. 54 (3), 84-87 (2010).

[23] Stepanov S.E., Tsyganok I.I., *Harmonic transforms of complete Riemannian manifolds*, Math. Notes, 100 (3), 465-471 (2016).

[24] **Catino G., Mastrolia P., Monticelli D., Rigoli M.**, *Analytical and geometric properties of generic Ricci Solitons*, Trans. Am. Math. Soc. **368** (11), 7533–7549 (2016).

[25] **Koboyashi S.**, *Transformation groups in differential geometry*, Springer-Verlag, Berlin and Heidelberg (1995).

[26] **Fernandez-Lopez M., Garcia-Rio E.**, *On gradient Ricci solitons with constant scalar curvature*, Proc. Amer. Math. Soc., **144**: 1, 369-378 (2016).

[27] **Stepanov S.E., Tsyganok I.I.**, *Vanishing theorems for harmonic mappings into nonnegativity curved smooth manifolds and their applications*, Manuscripta Math., **154** (1-2), 79-90 (2017).

[28] **Arnold V.I., Kresin B.A.**, *Topological methods in hydrodynamics*, Springer-Verlag, New York (1998).

[29] **Yau S.T.**, *Erratum: Some function-theoretic properties of complete Riemannian manifold and their applications to geometry*, (1976), 659-670, Indiana Univ. Math. J., **31** (4), 607 (1982).

[30] **O´Neil B.**, *Semi-Riemannian geometry with applications to Relativity*, Academic Press, San Diego (1983).

[31] **Davis W.R., Oliver D.R.**, *Matter field space times admitting symmetry mappings satisfying vanishing contraction of the Lie deformation of the Ricci tensor*, Ann. Inst. H. Poincare, Sec A (N.S.), **28** (2), 197-206 (1978).

[32] **Green L.H., Norris L.K., Oliver D.R., Davis W.R.**, *The Robertson-Walker metric and the symmetries belong to the family of contracted Ricci collineations*, General Relativity and Gravitation, **8** (9), 731-736 (1997).

[33] **Stepanov S.E., Shandra I.G.**, *New characteristics of infinitesimal isometry and nontrivial Ricci solitons*, Math. Notes **92** (3), 119-122 (2012).

[34] **Caminha A., Souza P., Camargo F.**, *Complete foliations of space forms by hypersurfaces*, Bull. Braz. Math. Soc., New Ser., **41** (3), 339–353 (2010).

[35] **Caminha A.**, *The geometry of closed conformal vector fields on Riemannian spaces*, Bull. Braz. Math. Soc., New Ser., **42** (2), 277–300 (2011).

[36] **Hamilton R. S.**, *Three-manifolds with positive Ricci curvature*, J. Diff. Geom., **17**: 2 (1982), 255-306.